\documentclass[reqno]{amsart}
\usepackage{amsmath,amsthm,amssymb}
\usepackage{url}
\newcommand{\be}{\begin{equation}}
\newcommand{\ee}{\end{equation}}
\newtheorem{lemma}{Lemma}
\newcommand{\di}{\displaystyle}
\begin{document}

\title[Heuristics for $\zeta '(-2k)$]{A heuristic derivation of  linear recurrence relations for $\zeta '(-2k)$ and $\zeta(2k+1)$}
\author{H. Gopalakrishna Gadiyar and R. Padma\\ ~~\\}
%\date{~~~}
\maketitle

\abstract 
We have gone back to old methods  found in the historical part of Hardy's Divergent Series well before the invention of the modern analytic continuation to use  formal manipulation of harmonic sums which produce some interesting formulae. These are linear recurrence relations for $\di{ \sum_{n=1}^\infty H_n n^k}$ which in turn yield  linear recurrence relations for $\zeta '(-k)$ and hence using the functional equation to a linear recurrence relation for  $\zeta '(2k)$ and $\zeta (2k+1)$. Questions of rigor have been postponed to a subsequent preprint. 

\endabstract

%\subjclass
%%%%%
%%%Mathematics Subject Classification 2010
%%%%%
%11A07, 11T71, 11Y16, 14G50, 68Q25, 94A60
%\endsubjclass

\section{Introduction}

For several years the authors have been studying the Euler sum 
\be
h(s) = \sum_{n=1}^\infty \frac{H_n}{n^s} \, .
\ee
where ${\di H_n=\sum_{k=1}^n \frac{1}{k}}$, the $n^{th}$ harmonic number. This has led to published works with K.N. Boyadzhiev \cite{boyadzhiev1:gadiyar:padma}, \cite{boyadzhiev2:gadiyar:padma} and B. Candelpergher \cite{candelpergher:gadiyar:padma} which was presented in the ICMZeta 2010 in Chennai. Also see \cite{candelpergher}.

We now turn to a curious formula in the paper with Boyadzhiev \cite{boyadzhiev1:gadiyar:padma} .

\begin{lemma} {\bf {(Fundamental Lemma)}} For $\Re(s)>1$ 

\be
\sum_{n=1}^\infty  \frac{H_n}{(n+1)^s} = \sum_{n=1}^\infty \frac{H_n}{n^s} -\zeta (s+1) \label{eq:Hn}
\ee
where $\zeta (s)$ is the Riemann zeta function.
\end{lemma}

\proof From the definition of $H_n$
\be
H_{n+1}=H_n +\frac{1}{n+1} \, .
\ee
Using this,
\begin{eqnarray}
\sum_{n=1}^\infty  \frac{H_n}{(n+1)^s}& =& \sum_{n=1}^\infty  \frac{H_{n+1}-\frac{1}{n+1}}{(n+1)^s} \\ \nonumber
~~ &=& \sum_{n=1}^\infty  \frac{H_{n+1}}{(n+1)^s} - \sum_{n=1}^\infty \frac{1}{(n+1)^{s+1}} \\ \nonumber
~~ &=& \left ( \sum_{n=1}^\infty  \frac{H_{n}}{n^s} -1 \right ) -\left ( \zeta(s+1) -1 \right ) \\ \nonumber
~~ &=& \sum_{n=1}^\infty  \frac{H_{n}}{n^s} - \zeta(s+1)
\end{eqnarray}

\endproof
Now we naively analytically continue the formula (\ref{eq:Hn}) and to our surprise we get new results. Replacing $s$ by $-s $ in( \ref{eq:Hn}) we get
\be
\sum_{n=1}^\infty  H_n(n+1)^s = \sum_{n=1}^\infty H_n n^s -\zeta (1-s) \label{eq:Hn1}
\ee

Using the functional equation of the Riemann zeta - function
\be
\zeta (1-s) = 2 (2\pi )^{-s} \Gamma (s) \cos \frac{\pi s}{2} \zeta (s) \, ,
\ee
(\ref{eq:Hn1}) becomes
\be
\sum_{n=1}^\infty  H_n(n+1)^s = \sum_{n=1}^\infty H_n n^s -2 (2\pi )^{-s} \Gamma (s) \cos \frac{\pi s}{2} \zeta (s)
\ee
For $s=2$ this becomes

\be
\sum_{n=1}^\infty  H_n(n+1)^2 - \sum_{n=1}^\infty H_n n^2 = \frac{1}{12}
\ee
where we have used  the formula $\di {\zeta (2) = \frac{\pi ^2}{6}}$. 
Expanding the left hand side and simplifying
\be
2~ \sum_{n=1}^\infty  n H_n +\sum_{n=1}^\infty H_n  = \frac{1}{12}
\ee
When $s=3$ we get
\be
3\sum_{n=1}^\infty  n^2 H_n +3 \sum_{n=1}^\infty n H_n + \sum_{n=1}^\infty  H_n= 0
\ee
as $\di {\cos \frac{3 \pi}{2}=0}$. Thus we get a recurrence relation for ${\di  \sum_{n=1}^\infty H_n n^k } $.

\section{Connection to $\zeta '(-k)$}
In this section we would like to give the relation between ${\di  \sum_{n=1}^\infty H_n n^k } $ and $\zeta '(-k)$. To do this we need an interpolation of the Bernoulli numbers. The Bernoulli numbers are defined as
\begin{equation}
\frac{z}{e^z-1}=\sum_{n=0}^\infty (-1)^n~B_n ~\frac{z^n}{n!} \label{eq:ber}
\end{equation}
With the definition given above  $\di {B_1(=\frac{1}{2})}$. Note that it differs by the conventional value of $\di {B_1(=-\frac{1}{2})}$ by a sign . The other values of Bernoulli numbers are not affected. These numbers can be got from the contour integral as follows.
\begin{lemma}
For $n \ge 1$,
\begin{equation}
\frac{1}{2\pi i} \int_C \frac{1}{z^{n}} \frac{-z}{e^{-z}-1}~ dz = \frac{B_{n-1}}{(n-1)!} \, .
\end{equation} 
where $C$ is the Hankel contour.
\end{lemma} 

Using the above integral we get an interpolation $B_s$ of the Bernoulli numbers. It is given by
\begin{equation}
-\frac{B_{s+1}}{\Gamma(s+2)} = \frac{1}{2\pi i} \int_C \frac{1}{z^{s+1}}~ \frac{e^z}{1-e^{z}}~dz \, . \label{eq:interpolation}
\end{equation}

The formulae connecting $B_s$ and its derivative $B_s ' $ with the Riemann zeta function are given below.
\begin{lemma}
\be
B_s = -s \zeta (1-s)
\ee

\be
B_s'=-\zeta(1-s)+s \zeta '(1-s)
\ee
\end{lemma}

In particular, 
\be
B'_{k}=-\zeta(1-k)+k \zeta '(1-k) \, , ~~k=1,2,\cdots \label{eq:B'k}
\ee
where  $\di{\zeta'(-k)=\sum_{n=1}^\infty n^{k} \log n }$ in an appropriate sense. 

The functional equation of $\zeta (s)$ gives the following relation between  $\zeta (2k+1)$ and $\zeta '(-2k)$

\be
\zeta (2k+1) = (-1)^k 2 \frac{(2\pi )^{2k}}{ (2k)! } \zeta ' (-2k)
\ee

Using this we get
\begin{lemma}
\be
\zeta (2k+1) = (-1)^k \frac{2\pi )^{2k+1}}{(2k+1)!} \frac{B'(2k+1)}{\pi},  
\ee
\end{lemma}

In \cite{gadiyar:padma} we proved the following identity.

\begin{lemma}
For $k \ge 1$,
\begin{equation} 
(-1)^{k-1}~k~\sum_{n=1}^\infty H_n n^{k-1}=B'_k+k B_{k-1}+\gamma B_k- \sum_{l+m=k} \left(\!\!\!
  \begin{array}{c}
    k \\
    l
  \end{array}
  \!\!\!\right) \frac{B_l}{l}B_m -B_k H_k \label{eq:matext}
\end{equation}
\end{lemma}
Susbstituting for $B_k '$ from (\ref{eq:B'k}) in terms of $\zeta '$, the above formula takes the form
\begin{equation} 
(-1)^{k-1}~k~\sum_{n=1}^\infty H_n n^{k-1}=-\zeta(1-k)+k \zeta '(1-k) +k B_{k-1}+\gamma B_k- \sum_{l+m=k} \left(\!\!\!
  \begin{array}{c}
    k \\
    l
  \end{array}
  \!\!\!\right) \frac{B_l}{l}B_m -B_k H_k \label{eq:matext1}
\end{equation}

\section{Comments}
\begin{itemize}
\item[1.] As it can be seen that the formulae will lead to a reduction in the number of unknown constants in the theory of the Riemann zeta function. If $\zeta '(0)$ is fed into the recurrence relations, all the $\zeta '(-k)$ for $k \ge 1$ will be obtained. 

\item[2.] We are trying to develop the theory of $\di{\sum_{n=1}^\infty H_n n^k}$ in strict analogy with Bernoulli numbers. In particular we want to get the analogue of 
\be
B_n = \sum_{k=0}^n   \left(\!\!\! \begin{array}{c}
    n \\
   k
  \end{array}\!\!\!\right) B_k \, ~{\rm for}  ~n \ge 2 
\ee

This is got from the identity 
\be
\frac{x e^x}{e^x-1} -\frac{x}{e^x-1} = x
\ee

\item[3.] An analogue of this for $\di {\sum_{n=1}^\infty H_n n^{k}}$ can be got by using the generating functions.
Let us start with the identity
\be
\frac{\log(1-e^{-x})}{1-e^{-x}} -\frac{e^{-x} \log(1-e^{-x})}{1-e^{-x}} = \log (1-e^{-x}) \label{eq:genfn}
\ee
Taking the Mellin transform of the above equation (\ref{eq:genfn}) 
\begin{eqnarray}
\int_0^\infty  x^{s-1} \frac{\log(1-e^{-x})}{1-e^{-x}} dx  &- & \int_0^\infty  x^{s-1}\frac{e^{-x} \log(1-e^{-x})}{1-e^{-x}} dx \nonumber \\
~~ &=& \int_0^\infty  x^{s-1} \log (1-e^{-x}) dx
\end{eqnarray}
and using the generating functions
\be
 \log (1-e^{-x}) = -\sum_{n=1}^\infty \frac{e^{-nx}}{n}
\ee
and
\be
\frac{\log(1-e^{-x})}{1-e^{-x}}=-\sum_{n=1}^\infty H_n~e^{-nx}
\ee
we get
\be
\sum_{n=1}^\infty  \frac{H_n}{(n+1)^s} = \sum_{n=1}^\infty \frac{H_n}{n^s} -\zeta (s+1)
\ee
which is our Fundmental Lemma (\ref{eq:Hn}).
\end{itemize}

\section{Conclusion} We paraphrase G. H. Hardy, see Section 1.3, \cite{hardy} (because we cannot imporove on his language):

 {\it The results of the formal calculations of  $\cdots $ are correct  wherever they can be checked: thus all of the formulae $\cdots $ are correct. It is natural to suppose that the other formulae will prove to be correct, and our transformations justifiable, if they are interpreted appropriately. We should then be able to regard the transformations as shorthand representations of more complex processes justifiable by the ordinary canons of analysis. It is plain that the first step towards such an interpretation must be some definition, or defintitions, of the `sum' of an infinite series, more widely applicable than the classical definition of Cauchy.

This remark is trivail now: it does not occur to a modern mathematician that a collection of mathematical symbols should have a `meaning' until one has been assigned to it by definition. It was not a triviality even to the greatest mathematicians of the eighteenth century. They had not the habit of definition: it was not natural to them to say, in so many words, `by X {\bf we mean} Y'. There are reservations to be made, to which we shall return in $\cdots $; but it is broadly ture to say that mathematicians before Cauchy asked not {\bf `How shall we define $1-1+1- \cdots $ ?'} but {\bf `What is $1 - 1+1 - \cdots ?$ ' }, and that this habit of mind led them into unnecessary perplexitites and controversies which were often really verbal.}

It is well known to experts in summability theory that different methods of summation assign slightly different values as the sum of a divergent series. All these points will have to be taken into account to get the value of $\zeta (2k+1)$.
\section{Acknowledgements} We would like to thank Professors Bernard Candelpergher and Khristo N. Boyadzhiev for many interesting discussions over email on the topic. We have also discussed the idea with Professor A. Sankaranaryanan, University of Hyderabad and Professor K. Srinivas, IMSc. All of them raised questions of rigor in the argument. We felt that the heuristic argument may be of wide interest and hence prepared this preprint. This preprint is meant to be a surprise to Prof. M. S. Rangachari, Ramanujan Institute for Advanced Study in Mathematics who kept alive research in divergent series in Chennai and constantly encouraged us to study the works of G. H. Hardy and S. Ramanujan.

\section{Dedication} We would like to dedicate this preprint to the memory of Prof. K. Ramachandra who encouraged study of simple ideas. We would also remember with gratitude our parents Hemalatha Gadiyar, Dr. H. G. Madhav Gadiyar,  R. Gowri and Prof. S. Ramanathan.

\noindent{Address of the authors}\\
48 Tulips Hermitage, Pennathur Road, Thirumalaikodi\\
Vellore 632055 Tamil Nadu INDIA\\
e-mail:gopigadiyar@gmail.com , padma.ramanathan@gmail.com
\end{document}